# Strong invariance principle for dependent random fields[*]

**Alexander Bulinski[1] and Alexey Shashkin[1]**

*Moscow State University*

**Abstract:** A strong invariance principle is established for random fields which satisfy dependence conditions more general than positive or negative association. We use the approach of Csörgő and Révész applied recently by Balan to associated random fields. The key step in our proof combines new moment and maximal inequalities, established by the authors for partial sums of multiindexed random variables, with the estimate of the convergence rate in the CLT for random fields under consideration.

## 1. Introduction and main results

Strong invariance principles are limit theorems concerning strong approximation for partial sums process of some random sequence or field by a (multiparameter) Wiener process. The first result of such type was obtained by Strassen [21] with the help of Skorokhod's embedding technique. Another powerful method, introduced by Csörgő and Révész [10], is based on quantile transforms. It was used by Kómlos, Major and Tusnady [15, 16] to achieve an unimprovable rate of convergence in the strong invariance principle for independent identically distributed random sequences. Berkes and Morrow [2] extended that method to mixing random fields.

In this paper, we study random fields with dependence condition proposed by Bulinski and Suquet [7] (in the case of a random sequence it was given by Doukhan and Louhichi [11]). Namely, let $X = \{X_j, j \in \mathbb{Z}^d\}$ be a real-valued random field on a probability space $(\Omega, \mathcal{F}, \mathsf{P})$ with $\mathsf{E}X_j^2 < \infty$ for any $j \in \mathbb{Z}^d$. We say that $X$ is *weakly dependent,* or $(BL, \theta)-$dependent, if there exists a sequence $\theta = (\theta_r)_{r \in \mathbb{N}}$ of positive numbers, $\theta_r \to 0$ as $r \to \infty$, such that for any pair of disjoint finite sets $I, J \subset \mathbb{Z}^d$ and any pair of bounded Lipschitz functions $f : \mathbb{R}^{|I|} \to \mathbb{R}$ and $g : \mathbb{R}^{|J|} \to \mathbb{R}$ one has

$$|cov(f(X_i, i \in I), g(X_j, j \in J))| \leq Lip(f)Lip(g)(|I| \wedge |J|)\theta_r. \qquad (1.1)$$

Here and below $|V|$ stands for the cardinality of a finite set $V$, $r = dist(I, J) = \min\{\|i - j\| : i \in I, j \in J\}$ with the norm $\|z\| = \max_{i=1,\ldots,d} |z_i|$, $z = (z_1, \ldots, z_d) \in \mathbb{Z}^d$, and, for $F : \mathbb{R}^n \to \mathbb{R}$,

$$Lip(F) = \sup_{x \neq y} \frac{|F(x) - F(y)|}{|x_1 - y_1| + \cdots + |x_n - y_n|}.$$

Note that one can apply (1.1) to unbounded Lipschitz functions $f$ and $g$ whenever $\mathsf{E}f^2(X_i, i \in I) < \infty$ and $\mathsf{E}g^2(X_j, j \in J) < \infty$.

---

[*]This work is partially supported by the RFBR grant 03-01-00724, by the grant 1758.2003.1 of Scientific Schools and by the INTAS grant 03-51-5018.

[1]Dept. of Mathematics and Mechanics, Moscow State University, Moscow, 119992, Russia, e-mail: bulinski@mech.math.msu.su e-mail: shashkin@mech.math.msu.su

The interest in studying model (1.1) is motivated by the following fact. There are a number of important stochastic models in mathematical statistics, reliability theory and statistical physics involving families of positively and negatively associated random variables (see [13, 14, 18] for the exact definitions and examples, for further references see, e.g., [7]). As shown by Bulinski and Shabanovich [8], a positively or negatively associated random field with finite second moments satisfies (1.1), provided that the Cox-Grimmett coefficient

$$\widehat{\theta}_r = \sup_{j \in \mathbb{Z}^d} \sum_{u \in \mathbb{Z}^d : \|u-j\| \geq r} |cov(X_u, X_j)|, \ r \in \mathbb{N},$$

is finite and $\widehat{\theta}_r \to 0$ when $r \to \infty$. In this case one can take $\theta = (\widehat{\theta}_r)_{r \in \mathbb{N}}$. There are also examples of (1.1) which are not induced by association, see [11, 20].

A strong invariance principle for associated random sequences whose Cox-Grimmett coefficient decreases exponentially was proved by Yu [22]. Recently Balan [1] extended this result to associated random fields. The principal goal of this paper is to extend the strong invariance principle to $(BL, \theta)$–dependent random fields. The new maximal inequality needed is given in Theorem 1.1.

For any finite $V \subset \mathbb{Z}^d$, we let $S(V) = \sum_{j \in V} X_j$. The sum over empty set is zero, as usual. We call a block a set $V = (a, b] := ((a_1, b_1] \times \cdots \times (a_d, b_d]) \cap \mathbb{Z}^d$ when $a, b \in \mathbb{Z}^d, a_1 < b_1, \ldots, a_d < b_d$. Given a block $V$, set $M(V) = \max_{W \subset V} |S(W)|$ where the supremum is over all the blocks $W$ contained in $V$.

Assume that
$$D_p := \sup_{j \in \mathbb{Z}^d} \mathsf{E}|X_j|^p < \infty \text{ for some } p > 2. \tag{1.2}$$

We will use condition (1.1) specialized to a sequence $\theta$ with a power or exponential rate of decreasing. Namely, either

$$\theta_r \leq c_0 r^{-\lambda}, \ r \in \mathbb{N}, \text{ for some } c_0 > 1 \text{ and } \lambda > 0, \tag{1.3}$$

or

$$\theta_r \leq c_0 e^{-\lambda r}, \ r \in \mathbb{N}, \text{ for some } c_0 > 1 \text{ and } \lambda > 0. \tag{1.4}$$

Introduce a function

$$\psi(x) = \begin{cases} (x-1)(x-2)^{-1}, & 2 < x \leq 4, \\ (3 - \sqrt{x})(\sqrt{x}+1)/2, & 4 < x \leq t_0^2, \\ ((x-1)\sqrt{(x-2)^2 - 3} - x^2 + 6x - 11)(3x-12)^{-1}, & x > t_0^2 \end{cases} \tag{1.5}$$

where $t_0 \approx 2.1413$ is the maximal root of the equation

$$t^3 + 2t^2 - 7t - 4 = 0.$$

Note that $\psi(x) \to 1$ as $x \to \infty$.

Now let us formulate the first of the main results of this paper.

**Theorem 1.1.** *Let $X$ be a centered $(BL, \theta)$–dependent random field satisfying (1.2) and (1.3) with $\lambda > d\psi(p)$ for $\psi(p)$ defined in (1.5). Then there exist $\delta > 0$ and $C > 1$ depending only on $d, p, D_p, c_0$ and $\lambda$ such that for any block $U \subset \mathbb{Z}^d$ one has*

$$\mathsf{E}|S(U)|^{2+\delta} \leq C|U|^{1+\delta/2}, \ \mathsf{E}M(U)^{2+\delta} \leq AC|U|^{1+\delta/2} \tag{1.6}$$

*where $A = 5^d(1 - 2^{\delta/(4+2\delta)})^{-d(2+\delta)}$.*

**Remark 1.** Moment and maximal inequalities for associated random fields were obtained in [4] and [6]. In the paper [19] similar inequalities were proved for weakly dependent random fields $X = \{X_j, j \in \mathbb{Z}^d\}$ under a stronger moment condition $\sup_{j \in \mathbb{Z}^d} \mathsf{E}|X_j|^{4+\delta} < \infty$, $\delta > 0$. An inequality for a $(BL, \theta)$−dependent field $X$ having only finite second moments when $\lambda > 3d$ in (1.3) was also established there, permitting to prove a weak invariance principle in the strictly stationary case. That result does not comprise ours.

**Remark 2.** The condition on the rate of decrease of $\theta_r$ determined by the function $\psi$ in (1.5) is implied by a simple condition $\lambda > d(p-1)/(p-2)$ because $\psi(p) \leq (p-1)/(p-2)$ for all $p > 2$.

Now suppose that
$$\sigma^2 := \sum_{j \in \mathbb{Z}^d} cov(X_0, X_j) \neq 0. \tag{1.7}$$

Note that (1.1) entails the convergence of series in (1.7) for a field $X$ with $\mathsf{E}X_j^2 < \infty$, $j \in \mathbb{Z}^d$. As is generally known, for a wide-sense stationary field $X$ one has
$$var(S_N) \sim \sigma^2[N] \text{ as } N \to \infty \tag{1.8}$$

where $N \in \mathbb{N}^d$, $[N] = N_1 \ldots N_d$, $S_N = S((0, N])$ and $N \to \infty$ means that $N_1 \to \infty, \ldots, N_d \to \infty$.

Following [2], for any $\tau > 0$, we introduce the set
$$G_\tau = \bigcap_{s=1}^{d} \left\{j \in \mathbb{N}^d : j_s \geq \left(\prod_{s' \neq s} j_{s'}\right)^\tau\right\}. \tag{1.9}$$

**Theorem 1.2.** *Suppose that $X$ is a wide-sense stationary $(BL, \theta)$−dependent centered random field satisfying (1.2), (1.4) and (1.7). Then one can redefine $X$, without changing its distribution, on a new probability space together with a $d$−parameter Wiener process $W = \{W_t, t \in [0, \infty)^d\}$, so that for some $\varepsilon > 0$ the following relation holds*
$$S_N - \sigma W_N = O([N]^{1/2-\varepsilon}) \text{ a.s.} \tag{1.10}$$

*as $N \to \infty$, $N \in G_\tau$ and $\tau > 0$.*

**Remark 3.** The value $\varepsilon$ in (1.10) depends on a field $X$. More precisely, $\varepsilon$ is determined by $\tau$, the covariance function of $X$ and parameters $d, p, D_p, c_0, \lambda$. Note that $G_\tau = \varnothing$ for $\tau > 1/(d-1)$, $d > 1$. One can easily obtain an analogue of Theorem 1.2 for wide-sense stationary weakly dependent stochastic process (i.e. for $d = 1$).

## 2. Proof of Theorem 1.1

We fix some $\delta \in (0, 1]$, $\delta < p - 2$. The exact value of $\delta$ will be specified later. Choose $A_\delta > 0$ (e.g., $A_\delta = 5$) to ensure that

$$(x+y)^2(1 + x + y)^\delta \leq x^{2+\delta} + y^{2+\delta} + A_\delta((1+x)^\delta y^2 + x^2(1+y)^\delta)$$

for any $x, y \geq 0$.

Let $h(n) = \min\{k \in \mathbb{Z}_+ : 2^k \geq n\}$, $n \in \mathbb{N}$. For any block $V \subset \mathbb{Z}^d$ having edges with lengths $l_1, \ldots, l_d$, we set $h(V) = h(l_1) + \cdots + h(l_d)$.

We will show that for some $C > 2(D_p \vee 1)$ and all blocks $U \subset \mathbb{Z}^d$

$$\mathsf{E}S^2(U)(1 + |S(U)|)^\delta \leq C|U|^{1+\delta/2}. \tag{2.1}$$

This is proved by induction on $h(U)$. For $h(U) = 0$ (i.e. when $|U| = 1$) inequality (2.1) is obviously true. Suppose now that (2.1) is verified for all $U$ such that $h(U) \leq h_0$. Consider a block $U$ having $h(U) = h_0 + 1$.

Let $L$ be any of the longest edges of $U$. Denote its length by $l(U)$. Draw a hyperplane orthogonal to $L$ dividing it into two intervals of lengths $[l(U)/2]$ and $l(U) - [l(U)/2]$, here $[\cdot]$ stands for integer part of a number. This hyperplane divides $U$ into two blocks $U_1$ and $U_2$ with $h(U_1), h(U_2) \leq h_0$.

**Lemma 2.1.** *There exists a value $\tau_0 = \tau_0(\delta) < 1$ such that, for any block $U \subset \mathbb{Z}^d$ with $|U| > 1$, one has*

$$|U_1|^{1+\delta/2} + |U_2|^{1+\delta/2} \leq \tau_0 |U|^{1+\delta/2}.$$

*Proof.* Straightforward. □

Observe that for the considered field $X$ condition (1.3) implies the bound

$$\mathsf{E} S^2(U) \leq (D_2 + c_0)|U| \qquad (2.2)$$

for any block $U \subset \mathbb{Z}^d$.

Set $Q_k = S(U_k)$, $k = 1, 2$. By induction hypothesis and Lemma 2.1,

$$\begin{aligned}
\mathsf{E} S^2(U)(1 + |S(U)|)^\delta &= \mathsf{E}(Q_1 + Q_2)^2(1 + |Q_1 + Q_2|)^\delta \leq C(|U_1|^{1+\delta/2} + |U_2|^{1+\delta/2}) \\
&\quad + A_\delta \mathsf{E}((1 + |Q_1|)^\delta Q_2^2 + (1 + |Q_2|)^\delta Q_1^2) \\
&\leq C\tau_0 |U|^{1+\delta/2} + A_\delta \mathsf{E}\left((1 + |Q_1|)^\delta Q_2^2 + (1 + |Q_2|)^\delta Q_1^2\right). \quad (2.3)
\end{aligned}$$

Our goal is to obtain upper bounds for $\mathsf{E}(1 + |Q_1|)^\delta Q_2^2$ and $\mathsf{E}(1 + |Q_2|)^\delta Q_1^2$. We proceed with the first estimate only, the second one being similar. To this end, let us take positive $\zeta < (1 - \tau_0)/(4A_\delta)$ and introduce a block $V$ letting

$$V = \{j \in U_2 : dist(\{j\}, U_1) \leq \zeta |U|^{1/d}\}.$$

Note that the induction hypothesis applies to $V$, since $V \subseteq U_2$. Using the Hölder inequality and (2.2), one shows that

$$\begin{aligned}
\mathsf{E}(1 + |Q_1|)^\delta Q_2^2 &\leq 2\mathsf{E}(1 + |Q_1|^\delta)S^2(V) \\
&\quad + 2\mathsf{E}(1 + |Q_1|)^\delta S^2(U_2 \setminus V) \\
&\leq 2(D_2 + c_0)|V| + 2(\mathsf{E}|Q_1|^{2+\delta})^{\delta/(2+\delta)}(\mathsf{E}|S(V)|^{2+\delta})^{2/(2+\delta)} \\
&\quad + 2\mathsf{E}(1 + |Q_1|)^\delta S^2(U_2 \setminus V) \\
&\leq 2(D_2 + c_0)|U| + 2C\zeta |U|^{1+\delta/2} \\
&\quad + 2\mathsf{E}(1 + |Q_1|)^\delta S^2(U_2 \setminus V). \quad (2.4)
\end{aligned}$$

Fix any indices $i, j \in U_2 \setminus V$ and assume first that $i \neq j$. Then $dist(\{j\}, \{i\} \cup U_1) = m > 0$. For any $y > 0$, we define the function $G_y$ by

$$G_y(t) = (|t| \wedge y) sign(t), \ t \in \mathbb{R}, \qquad (2.5)$$

and for some $y, z \geq 1$ introduce the random variables

$$Q_1^I = G_y(Q_1), \ Q_1^{II} = \left((1 + |Q_1|)^\delta - (1 + |Q_1^I|)^\delta\right)^{1/\delta},$$

$$X_i^I = G_z(X_i), \ X_i^{II} = X_i - X_i^I.$$

To simplify the notation we do not write $Q_{1,y}^I, Q_{1,y}^{II}, X_{i,z}^I$ and $X_{i,z}^{II}$. Obviously

$$|\mathsf{E}(1+|Q_1|)^\delta X_i X_j| \leq |\mathsf{E}(1+|Q_1^I|)^\delta X_i^I X_j| \\ + \mathsf{E}(1+|Q_1|)^\delta |X_i^{II} X_j| + \mathsf{E}|Q_1^{II}|^\delta |X_i^I X_j|. \quad (2.6)$$

Note that $\Phi(v,w) = (1+|G_y(v)|)^\delta G_z(w)$ is a bounded Lipschitz function with $Lip(\Phi) \leq 2y^\delta + z$. Since $X$ is a weakly dependent centered field, we can write

$$|\mathsf{E}(1+|Q_1^I|)^\delta X_i^I X_j| = |cov((1+|Q_1^I|)^\delta X_i^I, X_j)| \leq (2y^\delta + z)\theta_m. \quad (2.7)$$

Let $q$ be a positive number such that $1/q + \delta/(2+\delta) + 1/p = 1$, that is $q = p(2+\delta)/(2p-2-\delta) < p$. By the Hölder and Lyapunov inequalities,

$$\mathsf{E}\left(1+|Q_1|\right)^\delta |X_i^{II} X_j| \leq (\mathsf{E}|X_i^{II}|^q)^{1/q} D_p^{1/p}\left(1 + (\mathsf{E}|Q_1|^{2+\delta})^{\delta/(2+\delta)}\right)$$
$$\leq 2C^{\delta/(2+\delta)}|U|^{\delta/2}\left(\frac{D_p}{z^{p-q}}\right)^{1/q} D_p^{1/p}, \quad (2.8)$$

the last estimate being due to the induction hypothesis.

For $r \in (\delta, 2+\delta)$ to be specified later,

$$\mathsf{E}|Q_1^{II}|^\delta |X_i^I X_j| \leq z\mathsf{E}|Q_1^{II}|^\delta |X_j| \leq z\Big(|cov(|Q_1^{II}|^\delta, |X_j|)| + \mathsf{E}|X_j|\mathsf{E}|Q_1^{II}|^\delta\Big)$$
$$\leq z\delta y^{\delta-1}\theta_m + 2zy^{\delta-r} D_1 C^{r/(2+\delta)}|U|^{r/2}. \quad (2.9)$$

The last inequality follows from the induction hypothesis and the fact that the function $v \to ((1+|v|)^\delta - (1+y)^\delta)\mathbb{I}\{|v| \geq y\}$ is a Lipschitz one.

Now from (1.3) and (2.6)—(2.9), denoting $T = 2c_0(1 \vee D_p)$, we conclude that

$$|\mathsf{E}(1+|Q_1|)^\delta X_i X_j| \leq TC^{r/(2+\delta)}\Big((y^\delta+z)m^{-\lambda} + |U|^{\delta/2}z^{1-p/q} + zy^{\delta-r}|U|^{r/2}\Big). \quad (2.10)$$

Let $\beta, \gamma$ be positive parameters. Introduce $y = |U|^{1/2}m^{\beta\lambda}$, $z = m^{\gamma\lambda}$. Then in view of (2.10) we obtain

$$|\mathsf{E}(1+|Q_1|)^\delta X_i X_j| \leq TC^{r/(2+\delta)}\left(|U|^{\delta/2}\sum_{k=1}^{3} m^{-\lambda\nu_k} + m^{-\lambda\nu_4}\right)$$

where

$$\nu_1 = 1 - \delta\beta, \ \nu_2 = \gamma\left(\frac{p}{q}-1\right), \ \nu_3 = (r-\delta)\beta - \gamma, \ \nu_4 = 1 - \gamma.$$

Our next claim is the following elementary statement.

**Lemma 2.2.** *For each $d \in \mathbb{N}$, any block $U \subset \mathbb{Z}^d$, every $\nu > 0$ and arbitrary $i \in U$ one has*

$$\sum_{j \in U, j \neq i} \|i - j\|^{-\nu} \leq c(d,\nu) f(|U|,d,\nu) \quad (2.11)$$

*where $c(d,\nu) > 0$ and*

$$f(|U|,d,\nu) = \begin{cases} |U|^{1-\nu/d}, & 0 < \nu < d, \nu \notin \mathbb{N}, \\ (1+\ln|U|)|U|^{1-\nu/d}, & 0 < \nu \leq d, \nu \in \mathbb{N}, \\ 1, & \nu > d. \end{cases}$$

*Proof.* The case $d = 1$ is trivial. For $d \geq 2$, consider $U = (a,b] \cap \mathbb{Z}^d$. Without loss of generality we can assume that $l_1 \leq \cdots \leq l_d$ where $l_s = b_s - a_s$, $s = 1, \ldots, d$. It is easily seen that

$$\sum_{j \in U, j \neq i} \|i - j\|^{-\nu} \leq 3^d d \sum_{s=0}^{d-1} \prod_{1 \leq m \leq s} l_m \sum_{k=l_s+1}^{l_{s+1}} k^{d-s-1-\nu}$$

where $l_0 = 0$ and a product over an empty set is equal to 1. Using the well-known estimates for sums $\sum_{k=r_1}^{r_2} k^\gamma$ by means of corresponding integrals and the estimates

$$\Big(\prod_{1 \leq m \leq s} l_m\Big) l^{d-s-\nu} \leq \Big(\prod_{1 \leq m \leq d} l_m\Big)^{1-\nu/d}$$

for $l = l_s$ and $l = l_{s+1}$ we come to (2.11). The Lemma is complete. □

Now pick $r$ close enough to $2+\delta$ and $\beta, \gamma$ in such a way that $\lambda \nu_k > d$, $k = 1, 2, 3$. One can verify that this is possible if

$$\lambda > \lambda_1(d) = d\frac{(2+\delta)(2p - 4 - \delta)}{4(p - 2 - \delta)}. \quad (2.12)$$

Moreover, to have simultaneously $\lambda \nu_4 > (1 - \delta/2)d$ it suffices to require

$$\lambda > \lambda_2(d) = d\Big(\frac{2+\delta}{2(p-2-\delta)} + 1 - \frac{\delta}{2}\Big). \quad (2.13)$$

The condition imposed on $\lambda$ in Theorem 1.1 enables us to satisfy (2.12) and (2.13) taking $\delta$ small enough when $p \leq 4$, respectively $\delta = p - \sqrt{p} - 2$ when $4 < p \leq t_0^2$ and

$$\delta = \frac{2}{3}\Big(p - 2 - \sqrt{(p-2)^2 - 3}\Big)$$

otherwise.

For arbitrary $i \in U_2 \setminus V$, set

$$\overline{U_2} = \{j \in U_2 \setminus V : \|j - i\| \geq \zeta |U|^{1/d}\}, \widehat{U_2} = \{j \in U_2 \setminus V : \|j - i\| < \zeta |U|^{1/d}\}.$$

By Lemma 2.2, for any $i \in U_2 \setminus V$, we have

$$\Big|\sum_{j \neq i, j \in U_2 \setminus V} \mathsf{E}(1 + |Q_1|)^\delta X_i X_j\Big| \leq \Big|\sum_{j \in \overline{U_2}} \mathsf{E}(1 + |Q_1|)^\delta X_i X_j\Big|$$

$$+ \Big|\sum_{j \in \widehat{U_2}, j \neq i} \mathsf{E}(1 + |Q_1|)^\delta X_i X_j\Big|$$

$$\leq TC^{r/(2+\delta)}|U|^{\delta/2}\Big(4\zeta^{-\lambda \nu_0} + \sum_{k=1}^4 c(d, \lambda \nu_k)\Big), \quad (2.14)$$

here $\nu_0 = \max_{k=1,\ldots,4} \nu_k$.

Now we treat the case of $i = j \in U_2 \setminus V$. Obviously, one has $\delta p/(p-2) < 2 + \delta$. Therefore, by Hölder's inequality and induction hypothesis we infer that

$$\mathsf{E}(1 + |Q_1|)^\delta X_i^2 \leq (\mathsf{E}|X_i|^p)^{2/p}\Big(1 + \Big(\mathsf{E}|Q_1|^{\delta p/(p-2)}\Big)^{(p-2)/p}\Big)$$

$$\leq TC^{\delta/(2+\delta)}|U|^{\delta/2}. \quad (2.15)$$

From (2.14) and (2.15) one deduces that

$$\begin{aligned}
\mathsf{E}(1+|Q_1|)^\delta S^2(U_2 \setminus V) &\leq \Big|\sum_{\substack{i,j \in U_2 \setminus V \\ i \neq j}} \mathsf{E}(1+|Q_1|)^\delta X_i X_j\Big| + \sum_{i \in U_2 \setminus V} \mathsf{E}(1+|Q_1|)^\delta X_i^2 \\
&\leq |U| \max_{i \in U_2 \setminus V} \Big|\sum_{j \neq i, j \in U_2 \setminus V} \mathsf{E}(1+|Q_1|)^\delta X_i X_j\Big| \\
&\quad + TC^{\delta/(2+\delta)}|U|^{1+\delta/2} \\
&\leq MC^{r/(2+\delta)}|U|^{1+\delta/2}
\end{aligned} \qquad (2.16)$$

where $M = T(1 + 4\zeta^{-\lambda\nu_0} + \sum_{k=1}^4 c(d, \lambda\nu_k))$.

Employing (2.3), (2.4) and (2.16) we conclude that

$$\mathsf{E} S^2(U)(1+|S(U)|)^\delta \leq \Big(C\tau_0 + 4A_\delta(D_2 + c_0)4CA_\delta\zeta + 4C^{r/(2+\delta)}A_\delta M\Big)|U|^{1+\delta/2}.$$

The first assertion of the Theorem is now easily verified on account of (2.1) if $C$ is so large that

$$(1 - \tau_0 - 4A_\delta\zeta)C > 4A_\delta M C^{r/(2+\delta)} + 4A_\delta(D_2 + c_0).$$

The second assertion follows from the first one and the Moricz theorem [17]. □

## 3. Proof of Theorem 1.2

The proof adapts the approach of [1] and [2]. However, as the random field under consideration possesses a dependence property more general than association, we have to involve other results on normal approximation and partial sums behaviour. We also give simplified proofs of some steps.

Let $\alpha > \beta > 1$ be integers specified later. Introduce

$$n_0 = 0, \; n_l := \sum_{i=1}^l (i^\alpha + i^\beta), \; l \in \mathbb{N}. \qquad (3.1)$$

For $k \in \mathbb{N}^d$, put $k - 1 = (k_1 - 1, \ldots, k_d - 1)$ and $\mathsf{N}_k = (n_{k_1}, \ldots, n_{k_d})$. Set

$$\begin{aligned}
B_k &= (\mathsf{N}_{k-1}, \mathsf{N}_k], \; H_k = \prod_{s=1}^d (n_{k_s-1}, n_{k_s-1} + k_s^\alpha], \; I_k = B_k \setminus H_k, \\
u_k &= S(H_k), \sigma_k^2 = var(u_k), v_k = S(I_k), \tau_k^2 = var(v_k).
\end{aligned} \qquad (3.2)$$

We can redefine the random field $\{u_k\}$ on another probability space together with a random field $\{w_k, k \in \mathbb{N}^d\}$ of independent random variables such that $w_k \sim N(0, \tau_k^2)$ and the fields $\{u_k\}$ and $\{w_k\}$ are independent.

Further on we will denote by $C$ any positive factor which could depend only on $d, c_0, \lambda, p, D_p, \tau$ and the covariance function of the field $X$ except when specially mentioned. Occasionally when $C$ is a positive random variable, we write $C(\omega)$.[1] Now we pass to a number of lemmas. In their formulations, the requirements "(1.3) holds" or "(1.4) holds" mean that the field $X$ is $(BL, \theta)$-dependent with the decrease of $\theta$ as mentioned.

Due to weak stationarity and (1.7), we have $\sigma_k^2 > 0$ for any $k \in \mathbb{N}^d$. There exists $k_0 \in \mathbb{N}$ such that $\tau_k^2 > 0$ for all $k \in \mathbb{N}^d$ with $\min_{s=1,\ldots,d} k_s \geq k_0$. This is an immediate consequence of the following statement.

---
[1] A point $\omega$, in general, belongs to an extension of the initial probability space.

**Lemma 3.1.** *Let $X$ be a random field such that $D_2 < \infty$ and (1.4) holds. Then, for any finite union $V$ of disjoint blocks $V_q$, one has*

$$\sigma^2 - \frac{1}{|V|}var(S(V)) = O(l(V)^{-1/2}), \text{ as } l(V) \to \infty, \tag{3.3}$$

*where $l(V)$ is the minimal edge of all the blocks $V_q$.*

*Proof.* We have

$$\left|\sigma^2|V| - var(S(V))\right| = \left|\sum_{j \in V, k \notin V} cov(X_j, X_k)\right| = |\Sigma_1 + \Sigma_2|,$$

where $\Sigma_1$ is the sum over $k \notin V$ and all $j$'s such that $dist(\{j\}, \{\mathbb{Z}^d \setminus V\}) \geq l(V)^{1/2}$, and $\Sigma_2$ is taken over $k \notin V$ and the rest $j \in V$. By weak dependence $|\Sigma_1| \leq C|V|\exp\{-\lambda l(V)^{1/2}\}$. Furthermore, $|\Sigma_2| \leq 2^d C|V|l(V)^{-1/2}$ with $C = \sum_{k \in \mathbb{Z}^d} |cov(X_0, X_k)| < \infty$ and the Lemma is true. □

In what follows, we will consider only $k$ "large enough", that is having $\min_s k_s \geq k_0$. For such $k$ and $x \in \mathbb{R}$ let

$$\xi_k = \frac{u_k + w_k}{\sqrt{\sigma_k^2 + \tau_k^2}}, \ F_k(x) = \mathsf{P}(\xi_k \leq x). \tag{3.4}$$

Then $\xi_k$ has a density $f_k(x)$.

Analogously to [1] we introduce the random variables

$$\eta_k = \Phi^{-1}(F_k(\xi_k)), \ e_k = \sqrt{\sigma_k^2 + \tau_k^2}(\xi_k - \eta_k), \ k \in \mathbb{N}^d,$$

where $\Phi(x) = (2\pi)^{-1/2}\int_{-\infty}^{x} e^{-t^2/2}dt$.

Let $\rho = \tau/8$, $L$ be the set of all indices $i$ corresponding to the ("good") blocks $B_i \subset G_\rho$, and $H$ be the set of points in $\mathbb{N}^d$ which belong to some good block. For each point $N = (N_1, \ldots, N_d) \in H$, let $N^{(1)}, \ldots, N^{(d)}$ be the points defined as follows:

$$N^{(s)}_{s'} := N_{s'}, \ s' \neq s \text{ and } N^{(s)}_{s} := \min_{\substack{n \in H \\ n_{s'} = N_{s'}, s' \neq s}} n_s.$$

We consider also the sets

$$R_k = (M_k, \mathsf{N}_k] \text{ where } M_k = ((\mathsf{N}^{(1)}_k)_1, \ldots, (\mathsf{N}^{(d)}_k)_d), \ L_k = \{i : B_i \subset R_k\}.$$

Clearly,

$$\begin{aligned}S(R_k) &= \sum_{i \in L_k} e_i + \sum_{i \in L_k} \sqrt{|B_i|}\left(\sqrt{\frac{\sigma_i^2 + \tau_i^2}{|B_i|}} - \sigma\right)\eta_i \\ &+ \sum_{i \in L_k} \sigma\sqrt{|B_i|}\eta_i - \sum_{i \in L_k} w_i + \sum_{i \in L_k} v_i.\end{aligned} \tag{3.5}$$

If $V = (a, b] \cap \mathbb{Z}^d_+$ is a block in $\mathbb{Z}^d_+$ for some $a, b \in \mathbb{Z}^d_+$, and $\overline{V} = (a, b] \subset \mathbb{R}^d_+$, then set $W(V) := W(\overline{V})$. Here $W(\overline{V})$ is defined as usual (i.e. the signed sum of $W(t_l)$ where $t_l$, $l = 1, \ldots, 2^d$, are the vertices of $\overline{V}$).

**Lemma 3.2.** *If* (1.2), (1.3) *and* (1.7) *are satisfied, then for any* $k \in L$

$$\sup_{x \in \mathbb{R}} |F_k(x) - \Phi(x)| \leq C[k]^{-\alpha\mu}$$

*where* $\mu > 0$ *does not depend on* $k$.

*Proof.* If $X$ is a weakly dependent random field and we replace some of the variables $X_j$ with independent from $X$ and mutually independent random variables $Y_j$, then the random field $X' = \{X'_j, j \in \mathbb{N}^d\}$ obtained is again weakly dependent with the same sequence $\theta = (\theta_r)_{r \in \mathbb{N}}$ (cf. [12]). Moreover, note that if $X'_j = Y_j \sim N(0, \tau_k^2/|I_k|)$ for $j \in I_k$ and $X'_j = X_j$ for other $j \in \mathbb{Z}^d$ then for some $C > 0$,

$$\frac{1}{|B_k|} var\Big(\sum_{j \in B_k} X'_j\Big) \geq \frac{\sigma_k^2 + \tau_k^2}{2^d[k]^\alpha} \geq C \qquad (3.6)$$

for any $B_k \subset G_\rho$, $k \in \mathbb{N}^d$. This is true since the middle expression in (3.6) is always positive and in view of (1.8) and (3.3) has positive limit as $k \to \infty$, $k \in L$. Now the desired result follows from Theorem 4 of [9] by standard Bernstein's technique, analogously to Theorem 5 of that paper. □

**Lemma 3.3.** *Under conditions* (1.2), (1.3) *and* (1.7) *for any* $K \in (0, \sqrt{2\alpha\mu})$ *one has*

$$|\Phi^{-1}(F_k(x)) - x| \leq C[k]^{-\alpha\mu + K^2/2}$$

*if* $|x| \leq K\sqrt{\ln[k]}$. *Here* $C$ *may also depend on* $K$.

*Proof.* Follows the main lines of that given in [10] (Lemmas 2 and 3). □

**Lemma 3.4.** *If* (1.2), (1.4) *and* (1.7) *hold, then*

$$\mathsf{E} e_k^2 \leq C[k]^{\alpha - \varepsilon_0}$$

*for* $\varepsilon_0 = \alpha\mu\delta/(2 + 2\delta)$ *where* $\delta > 0$ *is the same as in Theorem 1.1 and* $\mu$ *is the same as in Lemma 3.2.*

*Proof.* In view of (2.2) our task is to show that $\mathsf{E}(\eta_k - \xi_k)^2 \leq C[k]^{-\varepsilon_0}$. To this end we take $K$ from Lemma 3.3 and write

$$\begin{aligned}
\mathsf{E}(\eta_k - \xi_k)^2 &= \mathsf{E}(\eta_k - \xi_k)^2 \mathbb{I}\{|\xi_k| \leq K\sqrt{\ln[k]}\} + \mathsf{E}(\eta_k - \xi_k)^2 \mathbb{I}\{|\xi_k| > K\sqrt{\ln[k]}\} \\
&\leq C[k]^{-\alpha\mu + K^2/2} + (\mathsf{E}|\eta_k - \xi_k|^{2+\delta})^{2/(2+\delta)} (\mathsf{P}(|\xi_k| > K\sqrt{\ln[k]}))^{\delta/(2+\delta)} \\
&\leq C[k]^{-\alpha\mu + K^2/2} + C(\mathsf{P}(|\eta_k| > K\sqrt{\ln[k]}))^{\delta/(2+\delta)} \\
&\quad + C|\mathsf{P}(|\xi_k| > K\sqrt{\ln[k]}) - \mathsf{P}(|\eta_k| > K\sqrt{\ln[k]})|^{\delta/(2+\delta)} \\
&\leq C([k]^{-\alpha\mu + K^2/2} + [k]^{-K^2\delta/2(2+\delta)} + [k]^{-\alpha\mu\delta/(2+\delta)})
\end{aligned}$$

by Lemmas 3.2 and 3.3 and Theorem 1.1. The optimization in $K$ yields the result. □

**Lemma 3.5.** *If* (1.2), (1.4), (1.7) *hold and* $\alpha - \beta \leq \varepsilon_0/4$, *then*

$$\sup_{x \in \mathbb{R}} |f_k(x) - f(x)| \leq C$$

*where* $f(x) = (2\pi)^{-1/2} e^{-x^2/2}$, $f_k$ *appeared after* (3.4) *and* $\varepsilon_0$ *is the same as in Lemma 3.4.*

*Proof.* Let $\varphi_k(t) = \mathsf{E}\exp\{it\xi_k\}$, $\widehat{\varphi_k}(t) = \mathsf{E}\exp\{itu_k/\sqrt{\sigma_k^2 + \tau_k^2}\}$, $\varphi(t) = e^{-t^2/2}$, where $t \in \mathbb{R}$, $i^2 = -1$. Note that

$$\varphi_k(t) = \widehat{\varphi_k}(t)\exp\left\{-\tau_k^2 t^2/2(\sigma_k^2 + \tau_k^2)\right\}.$$

By Lemma 3.4, for any $t \in \mathbb{R}$

$$|\varphi_k(t) - \varphi(t)| \leq \mathsf{E}|\exp\{it\xi_k\} - \exp\{it\eta_k\}| \leq |t|\mathsf{E}|\xi_k - \eta_k| \leq C|t|[k]^{-\varepsilon_0/2}.$$

Therefore for any $T > 0$ and any $x \in \mathbb{R}$, denoting $\nu_k = \tau_k^2/(\sigma_k^2 + \tau_k^2)$, one has

$$\begin{aligned}
|f_k(x) - f(x)| &\leq \frac{1}{2\pi}\int_\mathbb{R}|\varphi_k(t) - \varphi(t)|dt \\
&\leq \frac{1}{2\pi}\int_{|t|\leq T}|\varphi_k(t) - \varphi(t)|dt + \frac{1}{2\pi}\int_{|t|\geq T}\left(|\widehat{\varphi_k}(t)|e^{-t^2\nu_k/2} + \varphi(t)\right)dt \\
&\leq CT^2[k]^{-\varepsilon_0/2} + C\int_T^\infty e^{-t^2/2}dt + C\int_T^\infty e^{-\nu_k Tt/2}dt \\
&\leq C\left(T^2[k]^{-\varepsilon_0/2} + 1 + \frac{1}{T\nu_k}\exp\left\{-\frac{T^2\nu_k}{2}\right\}\right).
\end{aligned}$$

The lemma follows if we take $T = [k]^{\alpha-\beta}$. □

**Remark 4.** Clearly the condition $\alpha - \beta \leq \varepsilon_0/4$ is equivalent to

$$(\alpha/\beta)(1 - \mu\delta/8(1+\delta)) \leq 1. \tag{3.7}$$

Note that in Lemmas 3.4, 3.5 and Lemma 3.6 below one can replace condition (1.4) on $\lambda$ with that used in Theorem 1.1.

Two lemmas which follow can be proven analogously to Lemmas 3.6 and 3.9 in [1].

**Lemma 3.6.** *Suppose that (1.2), (1.4), (1.7) hold and $\alpha$ is so large that $\varepsilon_0 > 2$. Then there exists $\varepsilon_1 > 0$ such that*

$$\sum_{i \in L_k}|e_i| \leq C(\omega)[\mathsf{N}_k]^{1/2-\varepsilon_1} \ a.s.$$

**Lemma 3.7.** *Assume that $D_2 < \infty$. If (1.3) holds and $\alpha - \beta > 6/\rho$ then, for some $\varepsilon_2 > 0$,*

$$\sum_{i \in L_k}(|v_i| + |w_i|) \leq C(\omega)[\mathsf{N}_k]^{1/2-\varepsilon_2} \ a.s.$$

**Lemma 3.8.** *Suppose that $D_2 < \infty$. If (1.4), (1.7) hold and $\beta > 6/\rho$ then, for some $\varepsilon_3 > 0$,*

$$\sum_{i \in L_k}\sqrt{|B_i|}\left(\sigma - \sqrt{\frac{\sigma_i^2 + \tau_i^2}{|B_i|}}\right)|\eta_i| \leq C(\omega)[\mathsf{N}_k]^{1/2-\varepsilon_3} \ a.s. \tag{3.8}$$

*Proof.* With the help of an inequality

$$\left|\sigma - \sqrt{\frac{\sigma_i^2 + \tau_i^2}{|B_i|}}\right| \leq \sigma^{-1}\left|\sigma^2 - \frac{1}{|B_i|}\mathrm{var}(u_i + w_i)\right| \leq C[i]^{-\beta\rho/4}$$

ensuing from Lemma 3.1, we come to (3.8) applying arguments analogous to those proving Lemma 3.8 in [1]. □

**Lemma 3.9.** *Suppose that $X = \{X_j, j \in \mathbb{Z}^d\}$ is a weakly dependent random field and $Y = \{Y_j, j \in \mathbb{Z}^d\}$ is a field consisting of independent random variables and independent of $X$. Let $I, J \subset \mathbb{Z}^d$ be disjoint finite sets and $f : \mathbb{R}^{|I|} \to \mathbb{R}$, $g : \mathbb{R}^{|J|} \to \mathbb{R}$ be bounded Lipschitz functions. Then*

$$|cov(f(X_I + Y_I), g(X_J + Y_J))| \leq Lip(f)Lip(g)(|I| \wedge |J|)\theta_r \qquad (3.9)$$

*for $r = dist(I, J)$, that is, such addition of $Y$ does not alter the property* (1.1).

*Proof.* By a smoothing procedure we can reduce the general case to that of the random vector $(Y_I, Y_J)$ with a density $q(t^1, t^2) = q_I(t^1)q_J(t^2)$, here $t^1 \in \mathbb{R}^{|I|}, t^2 \in \mathbb{R}^{|J|}$. Evidently (3.9) is true for a field $Y$ consisting of some constants. Thus by independence hypothesis, the Fubini theorem and (1.1) we have

$$|cov(f(X_I + Y_I), g(X_J + Y_J))|$$
$$= \left| \iint cov(f(X_I + y_I), g(X_J + y_J))q(y_I, y_J)dy_I dy_J \right|$$
$$\leq Lip(f)Lip(g)(|I| \wedge |J|)\theta_r \iint q(y_I, y_J)dy_I dy_J$$
$$= Lip(f)Lip(g)(|I| \wedge |J|)\theta_r$$

where the double integral is taken over $\mathbb{R}^{|I|} \times \mathbb{R}^{|J|}$. $\square$

**Lemma 3.10 (Lemma 4.3 in [1]).** *There exists a bijection $\psi : \mathbb{Z}_+ \to L$ with the following properties:*

$$l < m \Rightarrow \exists s^* = s^*(l, m) \in \{1, \ldots, d\} \text{ such that } \psi(l)_{s^*} \leq \psi(m)_{s^*},$$

$$\exists m_0 \in \mathbb{Z}_+ \text{ such that } m \leq [\psi(m)]^{\gamma_0} \; \forall m \geq m_0,$$

*for any $\gamma_0 > (1 + 1/\rho)(1 - 1/d)$.*

Set $Y_m = \eta_{\psi(m)}$, $m \in \mathbb{Z}_+$.

**Lemma 3.11.** *If* (1.2), (1.4), (1.7) *and* (3.7) *hold, then for every $m \in \mathbb{N}, m > 1$, and all $t = (t_1, \ldots, t_m) \in \mathbb{R}^m$ one has*

$$\left| \mathsf{E} \exp\left\{i \sum_{l=1}^{m} t_l Y_l\right\} - \mathsf{E} \exp\left\{i \sum_{l=1}^{m-1} t_l Y_l\right\} \mathsf{E} \exp\{it_m Y_m\} \right|$$
$$\leq Cm\Delta^\gamma, \; i^2 = -1, \qquad (3.10)$$

*where*

$$\Delta = \frac{[\psi(m)]^\alpha}{m} \theta_r \sum_{l=1}^{m} t_l^2, \; r = dist\left(H_{\psi(m)}, \cup_{l=1}^{m-1} H_{\psi(l)}\right), \; \gamma = \begin{cases} 1, & \Delta > 1, \\ 1/3, & \Delta \leq 1 \end{cases} \qquad (3.11)$$

*and the cubes $H_k$ were defined in* (3.2).

*Proof.* Let $M > 0$ be a number to be specified later and $G_M(t)$ be the function defined in (2.5). We set $Y_{j,M} = G_M(\eta_{\psi(j)})$, $j \in \mathbb{Z}_+$. Note that

$$|\exp\{itY_{j,M}\} - \exp\{itY_j\}| \leq 2I\{|Y_j| > M\}, \; t \in \mathbb{R}.$$

Therefore, we have

$$\left|cov\left(\exp\Big\{i\sum_{l=1}^{m-1}t_lY_l\Big\},\exp\{-it_mY_m\}\right)\right|$$
$$\leq \left|cov\left(\exp\Big\{i\sum_{l=1}^{m-1}t_lY_{l,M}\Big\},\exp\{-it_mY_{m,M}\}\right)\right|$$
$$+ 4\sum_{l=1}^{m}\mathsf{P}\{|Y_l|>M\}. \qquad (3.12)$$

Every summand in (3.12) except for the first one is not greater than $Ce^{-M^2/2}$. To estimate the first summand in the right hand side of (3.12), we notice that, for any $k \in \mathbb{N}^d$, the random variable $\eta_{k,M}$ is a Lipschitz function of $\xi_k$. Indeed, $Lip(F_k) \leq C$ by Lemma 3.5 and $\eta_{k,M} = h_M(F_k(\xi_k))$ where

$$h_M(x) = \big(|\Phi^{-1}(x)| \wedge M\big)\, sign(\Phi^{-1}(x)),\ x \in \mathbb{R}.$$

Clearly $Lip(h_M) \leq \sqrt{2\pi}e^{M^2/2}$. By Lemma 3.9, one can estimate the covariance in the same way as if the normal variables $w_k, k \in \mathbb{N}^d$ were constants, and with the help of (3.6) we obtain

$$\left|cov\left(\exp\Big\{i\sum_{l=1}^{m-1}t_lY_{l,M}\Big\},\exp\{-it_mY_{m,M}\}\right)\right| \leq C\theta_r e^{M^2}[\psi(m)]^{\alpha}\sum_{l=1}^{m}t_l^2. \qquad (3.13)$$

Thus, from (3.12)—(3.13) we see that

$$\left|cov\left(\exp\Big\{i\sum_{l=1}^{m-1}t_lY_l\Big\},\exp\{-it_mY_m\}\right)\right| \leq C\Big(me^{-M^2/2} + e^{M^2}[\psi(m)]^{\alpha}\theta_r\sum_{l=1}^{m}t_l^2\Big).$$

The result follows now by optimization in $M$. $\square$

**Lemma 3.12.** *Suppose that (1.2), (1.4), (1.7), (3.7) hold and $\beta > 2\gamma_0/\rho$ where $\gamma_0$ appears in the formulation of Lemma 3.10. Then we can redefine the random field $X$, without changing its distribution, on a new probability space together with a $d$−parameter Wiener process $W = \{W_t, t \in [0,\infty)^d\}$, such that for some $\varepsilon_4 > 0$*

$$\sum_{i\in L_k}\sigma\sqrt{|B_i|}\left|\eta_i - \frac{W(B_i)}{\sqrt{|B_i|}}\right| \leq C(\omega)[\mathsf{N}_k]^{1/2-\varepsilon_4}\ a.s. \qquad (3.14)$$

*Proof.* Arguing as in the proof of Lemma 4.4 of [1], we see that by Berkes-Philipp strong approximation result ([3]) it is enough to establish the following fact. There exist sequences $\varkappa_m > 0$ and $z_m > 10^4 m^2$ ($m \in \mathbb{N}$, $m > 1$) such that

(1) for any $m \in \mathbb{N}, m > 1$, and all $t = (t_1, \ldots, t_m) \in \mathbb{R}^m$ with $\sum_{l=1}^{m} t_l^2 \leq z_m^2$ one has

$$\left|\mathsf{E}\exp\Big\{i\sum_{l=1}^{m}t_lY_l\Big\} - \mathsf{E}\exp\Big\{i\sum_{l=1}^{m-1}t_lY_l\Big\}\mathsf{E}\exp\{it_mY_m\}\right| \leq \varkappa_m,$$

(2) $z_m^{-1/4}\ln z_m + \exp\{-3z_m^{1/2}/16\}m^{1/2}z_m^{1/4} + \varkappa_m^{1/2}z_m^{m+1/4} = O(m^{-2}), m \to \infty$.

We take $z_m = 10^4 m^q, q > 8$. Then it suffices to prove that one can take $\varkappa_m = O(\exp\{-Am^R\}), m \to \infty$, for some $A > 0$, $R > 1$.

The distance $r$ between $H_m$ and any of the blocks $H_1, \ldots, H_{m-1}$ is by construction not less than $\min_{s=1,\ldots,d}(\psi(m)_s - 1)^\beta$. But since $\psi(m) \in L$, by Lemma 3.10 for $m > m_0$ we have

$$r \geq C \min_{s=1,\ldots,d} \psi(m)_s^\beta \geq C[\psi(m)]^{\rho\beta/2} \geq Cm^{\rho\beta/2\gamma_0}. \tag{3.15}$$

From Lemma 3.10 one can also easily see that for large enough $m$

$$[\psi(m)] \leq C \min_{s=1,\ldots,d} \psi(m)_s^{2/\rho} \leq Cm^{2/\rho}. \tag{3.16}$$

Obviously, $\Delta_1 := Cm^{2q+2\alpha/\rho}\exp\{-C\lambda m^{\rho\beta/2\gamma_0}\} < 1$ for all $m$ large enough. Therefore, for such $m$, by Lemma 3.11, (1.4), (3.15) and (3.16), one has

$$\left| \mathsf{E}\exp\Big\{i\sum_{l=1}^m t_l Y_l\Big\} - \mathsf{E}\exp\Big\{i\sum_{l=1}^{m-1} t_l Y_l\Big\}\mathsf{E}\exp\{it_m Y_m\}\right| \leq Cm\Delta^\gamma \leq Cm\Delta_1^{1/3}$$

where $\Delta$ and $\gamma$ are defined in (3.11). Thus one can take $A = 1$, $R \in (1, \rho\beta/2\gamma_0)$. The lemma is proved. $\square$

Now we estimate the terms $S((0, \mathsf{N}_k] \setminus R_k)$, $\sigma W((0, \mathsf{N}_k] \setminus R_k)$, $S_N - S((0, \mathsf{N}_k])$, $\sigma W_N - \sigma W((0, \mathsf{N}_k])$ when $N \in G_\tau$ and $\mathsf{N}_k < N \leq \mathsf{N}_{k+1}$. Here the relation $a < b$ ($a \leq b$) for $a, b \in \mathbb{Z}^d$ is defined in the usual way. All the terms involving the Wiener process can be considered as sums of independent identically distributed normal random variables, therefore forming a weakly dependent field, we will proceed only with the partial sums generated by $X$. Obviously one can write

$$|S((0, \mathsf{N}_k] \setminus R_k)| \leq \sum_{s=1}^d 2^{d-s} \mathsf{M}_s(\mathsf{N}_k), \quad (0, N] \setminus (0, \mathsf{N}_k] = \cup_J I_k^{(J)}$$

where $\mathsf{M}_s(N) = \max_{n \leq N^{(s)}} |S_n|$, $I_k^{(J)} = \prod_{s\in J}(n_{k_s}, N_s] \times \prod_{s\notin J}(0, n_{k_s}]$ and the union is taken over all non-empty subsets $J$ of $\{1, \ldots, d\}$. Furthermore,

$$\max_{\mathsf{N}_k < N \leq \mathsf{N}_{k+1}} |S_N - S((0, \mathsf{N}_k])| \leq \sum_J M_k^{(J)} \text{ where } M_k^{(J)} = \sup |S(I_k^{(J)})|,$$

here the supremum is taken over all $N$ such that $n_{k_s} < N_s \leq n_{k_s+1}$, $s \in J$.

**Lemma 3.13.** *Suppose that conditions of Theorem 1.1 are satisfied. Then there exists $\delta > 0$ such that for any $x > 0$ and any block $V$,*

$$\mathsf{P}(M(V) \geq x\sqrt{|V|}) \leq Cx^{-2-\delta}.$$

*Proof.* Follows from the second assertion of Theorem 1.1 and the Markov inequality. $\square$

**Lemma 3.14.** *Let conditions of Theorem 1.1 hold. If (1.7) is true then there exists $\gamma_1 > 0$ such that for any block $V = (m, m+n]$ with $n \in G_\rho$ and $m \in \mathbb{Z}_+^d$*

$$\mathsf{P}(M(V) \geq |V|^{1/2}(\ln|V|)^{d+1}) \leq C|V|^{-\gamma_1},$$

*where $C$ does not depend on $m$ and $n$.*

*Proof.* The proof is the same as that for the second inequality of Lemma 7 in [2]; the needed Berry-Esseen type estimate for $\mathsf{P}(S(V) \geq x\sqrt{|V|})$, $x > 0$, can be obtained from the results mentioned in the proof of Lemma 3.2. □

The next two lemmas are proved analogously to Lemmas 6 and 9 established in [2] for mixing random fields.

**Lemma 3.15.** *Assume that conditions of Theorem 1.1 are satisfied and $\alpha > 8/(3\tau) - 1$. Then*

$$\max_{s=1,\ldots,d} \mathsf{M}_s(\mathsf{N}_k) \leq C(\omega)[\mathsf{N}_k]^{1/2-\varepsilon_5} \, a.s.$$

*for some $\varepsilon_5 > 0$ and every $\mathsf{N}_k \in G_\tau$.*

**Lemma 3.16.** *Suppose that conditions of Theorem 1.1 hold and (1.7) is true. Let $\gamma_1$ be the constant given by Lemma 3.14. If $\alpha > 2/\gamma_1$ then one has*

$$\max_J M_k^{(J)} \leq C(\omega)[\mathsf{N}_k]^{1/2-\varepsilon_6} \, a.s.$$

*for some $\varepsilon_6 > 0$ and every $\mathsf{N}_k \in G_\rho$.*

Now, if we take $\gamma_0, \alpha, \beta$ to satisfy conditions

$$\gamma_0 > \left(1 + \frac{1}{\rho}\right)\left(1 - \frac{1}{d}\right), \frac{\alpha}{\beta}\left(1 - \frac{\mu\delta}{8(1+\delta)}\right) < 1, \beta > \frac{6}{\rho},$$

$$\alpha - \beta > \frac{6}{\rho}, \beta > \frac{2\gamma_0}{\rho}, \alpha > \frac{8}{3\tau} - 1, \alpha\gamma_1 > 2,$$

then Theorem 1.2 follows from (3.5), Lemmas 3.6, 3.7, 3.8, 3.12, 3.15 and 3.16 with $\varepsilon = \min_{i=1,\ldots,6} \varepsilon_i$. The proof of Theorem 1.2 is completed. □

## 4. The law of the iterated logarithm

In order to provide applications of Theorem 1.2, we state now the law of the iterated logarithm for weakly dependent random fields. This is the first result of such type for the fields with dependence condition (1.1); it generalizes the laws of the iterated logarithm known for positively and negatively associated random variables, see, e.g., [5].

For $x > 0$, set $\mathrm{Log}\, x = \log(x \vee e)$.

**Theorem 4.1.** *Suppose that $X$ is a random field satisfying all the conditions of Theorem 1.2. Then, for any $\tau \in (0, 1/(d-1))$, one has, almost surely,*

$$\limsup \frac{S_N}{\sqrt{2d\sigma^2[N]\mathrm{LogLog}[N]}} = 1 \; and \; \liminf \frac{S_N}{\sqrt{2d\sigma^2[N]\mathrm{LogLog}[N]}} = -1,$$

*as $N \to \infty$, $N \in G_\tau$.*

*Proof.* We have $\limsup W_N/\sqrt{2d\sigma^2[N]\mathrm{LogLog}[N]} \leq 1$, due to the LIL for $d$−parameter Wiener process [23]. The fact that this upper limit is not less than 1 as $N \to \infty$, $N \in G_\tau$, can be proved in the same way as the lower bound in classical law of the iterated logarithm for the Wiener process. Thus, Theorem follows from Theorem 1.2. □

**Final remarks**

It should be noted that initially the authors placed on top of the first page the following dedication:

> To Professor M. S. Keane on the occasion of his anniversary.

To unify the style of the volume such inscriptions were omitted. We hope that Mike Keane will read this paper to discover that phrase, as well as our congratulations and best wishes.

**Acknowledgments**

A.Bulinski is grateful to the Department of Probability Theory and Stochastic Models of the University Paris-VI for hospitality.